# Comments to Neutrosophy

Carlos Gershenson
School of Cognitive and Computer Sciences
University of Sussex
Brighton, BN1 9QN, U. K.
C.Gershenson@sussex.ac.uk

*Abstract*

*Any system based on axioms is incomplete because the axioms cannot be proven from the system, just believed. But one system can be less-incomplete than other. Neutrosophy is less-incomplete than many other systems because it contains them. But this does not mean that it is finished, and it can always be improved. The comments presented here are an attempt to make Neutrosophy even less-incomplete. I argue that less-incomplete ideas are more useful, since we cannot perceive truth or falsity or indeterminacy independently of a context, and are therefore relative. Absolute being and relative being are defined. Also the "silly theorem problem" is posed, and its partial solution described. The issues arising from the incompleteness of our contexts are presented. We also note the relativity and dependance of logic to a context. We propose "metacontextuality" as a paradigm for containing as many contexts as we can, in order to be less-incomplete and discuss some possible consequences.*

## 1. Introduction

Upset because none of the logics I knew were able to handle contradictions and paradoxes, I created my own, which was able to do so (Gershenson 1998; 1999). Later I would find out that my Multidimensional Logic (MDL) fitted in the category of paraconsistent logics (Priest and Tanaka, 1996). Basically, I rejected the axiom of no contradiction, and instead of having truth values of propositions, I defined truth vectors, where each element was independent of the other. In this way, we can have a vector representing something that is true and false at the same time, or nor true nor false at the same time. MDL gave interesting results, and the new perspective influenced me to propose ideas in different branches of Philosophy, Complexity, and Computer Science (e.g. Gershenson, 1999; 2002).

In the autumn of 2001, I became aware of Neutrosophy (Smarandache, 1995). The fact that it reached similar results than the ones I had makes me think that there is a lower probability that we are completely mistaken. Neutrosophy covers a much wider area than my ideas, but nevertheless I believe it could be enriched by them, since there are non-overlapping parts. In this work I will expose how some of these ideas might help in making Neutrosophy less-incomplete.

The contents of this paper were mainly sparse comments made to Florentin Smarandache, and while giving them the shape of an article, the continuity of ideas was not an entire success. My apologies to the sensible reader.

## 2. To be *and* not to be...

When people speak about the being, since it is one of the most general things you can speak about, they often seem to speak about different things. I define two types of being: absolute and relative. The **absolute being** (a-being) is independent of the observer, infinite. The **relative being** (re-being) is dependant of the observer, therefore finite, and different in each individual. The re-being can approach as much *as we want to* towards the a-being, but it can never comprehend it. Objects a-are. Concepts, representations, and ideas re-are (Objects do not depend on the representations we have of them).

We can only suppose about what things a-are, we cannot be absolutely sure, we can only speculate, because they a-are infinite and we are not. We cannot say that something is absolutely true or false. We can only assert things in a relative way. We could assign truth values or vectors to them, but these would be relative to our **context**.

The **being** would be the conjunction of re-being and a-being. This is, without the distinction we are drawing between them. But isn't it confusing to speak about something which is absolute and relative at the same time? Yes, precisely. But that is what we do every day.

We can say that every effect a-has cause(s), but it does not mean that we can perceive them, so some effects do not re-have causes.

## 3. Reason, beliefs, and experience

Reason cannot prove itself. It would be as if a theorem would like to prove the axioms it is based upon. Reason, as all systems based in axioms, is incomplete (in the sense of Gödel (1931) and Turing (1936)). Beliefs are the axioms of reason and thought. Reason cannot prove the beliefs it is based upon. But how do these beliefs arise, then? Through experience. But experience needs of previous beliefs and reason to be assimilated, and reason also needs of experience to be formed, as beliefs need of reason as well. Beliefs, reason, and experience, are based upon each other. Which came first? All of them. We cannot have one without the others. **Contexts** are dynamic, and formed upon beliefs, reason, and experience. It is there where the re-being lies. Since the re-being is dependant of our context, it is also dependant of our beliefs, reasonings, and experiences. Contexts are dynamic because they are changed constantly as we have new experiences, change our beliefs, and our ways of reasoning.

For example, we cannot say if capitalism a-is good or bad (independent of a context). It re-is good for the people who get a profit from it, re-is neutral for the people who think that are not affected by it, and bad for people which suffer from it. There re-is a god for people who believe so, and there re-is not for people who do not believe so, but we cannot say if there a-is a god. We can speculate as much as we want to about everything which a-is, but we will never contain it, therefore we can only have *an idea* of it. It is not only that "we know without knowing" (Smarandache, 1995), but less-incompletely: we re-know, but a-know not.

Anything we want to assert has an implicit **context**, for which what we want to assert is consistent. But *every idea is valid in the context it was created*. Ideas cannot a-be right or wrong. They re-are right or wrong according to a context. And since in order to be created they need to be consistent with their context, they re-are always right according to their context. Then which idea is more valuable? We should try to see which context is more valuable. Well, we cannot say which a-is more valuable, we can only see that if a context contains others, the ideas created in it will be valid also in the other contexts, and not necessarily vice versa. We do not know which one a-is better, but the context which contains the others will be **less-incomplete**.

A problem arises on the horizon...

## 4. The silly theorem problem

For any silly theorem $T_s$, we can find at least one set of axioms such that $T_s$ is consistent with the system defined by the axioms. How do we know a theorem is not silly? Or how do we know if the axioms are not silly?

We can extend this same problem to contexts: there can re-be any silly idea $I_s$ so that there is at least one context for which it is consistent[1]. Empiricism comes to the rescue. Well, not only experience, but **evolution** more or less helps us get out of the problem. This is, contexts which support silly ideas, even when they re-are not silly in their own context, will not be able to contend with the contexts which are able to describe more closely (less-incompletely) what things a-are. The evolution of contexts consists in making them less-incomplete. In formal systems a similar evolutionary (tending to pragmatic) criterium could be seen: the axioms which support silly theorems are not useful, and that is why we then call them silly. If the axioms and the theorems derived from them are useful in a context, we do not consider them silly. But this is not a complete solution for the silly theorem problem. There is no complete solution for the silly theorem problem (Or by saying this *that* is the solution?).

But silly ideas (for our context) are useful as well, because our ideas become stronger when they neutralize them, and become less-incomplete if they contain them. There would not be smart ideas without silly ones.

We can take advantage of the silly theorem problem, in order to define axioms after the theorems, or contexts after ideas. For example, we might want such sets that A contains B, and B contains A, but without A and B being equal. We would just need to define the proper axioms. If this new set theory will be useful or silly is no reason for stopping, because we will not know if it a-is silly or not, and we will know if it re-is silly only after we create it.

---

[1] But the silliness is also relative to the context of the one who judges the idea...

## 5. Incomplete Language

Not only silliness, but all adjectives can only be <used|applied> in a relative way, dependant of a context. Language is relative as well. How can we speak about absolute being, then? We can and we cannot. We speak about it, but in that moment it a-is relative. For us, it is and it is not-incomplete. But that we cannot completely speak about it, it is not a reason to stop speaking about it (as Wittgenstein (1918) would early suggest in his Tractatus Logicus Philosophicus), because we can incompletely represent its completeness... As Wittgenstein himself (but not most of his followers...) realized, following the ideas in the Tractatus, we would not be able to speak about anything... (languages are incomplete). Language is used inside a context. Depending of this context the language will be different.

## 6. Incompleteness

If we cannot create complete systems, we should try to make them as less-incomplete as possible. Since our systems are incomplete because they are finite, there is no way of measuring the completeness of a system. It is like asking: how infinite is x, when x is a finite number? We could also have a huge system, but if it is not related to the a-being of something, it is not so useful as a small system which describes closely a part of what something a-is. So, we could say that a system **re-is** less-incomplete as it approaches more the a-being.

But, for example, what is the a-being of the number four? We define numbers, so we determine what they are. But this definition arises on our generalizations of what things a-are. There a-is no number four, but we could say that there the number four a-is in things.

The a-being is far from materialism. Materialism re-is, and we cannot say for sure if things a-are only materialistic or not, because we do not know what matter a-is.

So, getting back to the incompleteness of systems, they can approach the a-being as much as we want to, but we <do not|will never> know what the *a-being a-is*. We could measure the incompleteness of a system in a relative way, but I am not interested in defining such a method. But we can see that a system will be less-incomplete if it contains others, even the ones that (re-)"are wrong". This is because a system which tries not only to explain why things a-are, but also why other people thought it was something else will have a wider perspective. Because people a-are not mistaken (nor right). There cannot be errors inside their context (people do not create systems silly/mistaken for themselves). So, by containing as many contexts as possible, we are also containing as many systems as possible. The great attempts in science of unifying theories could be seen with these eyes. But for containing as many contexts as we can, we are required to leave all hope of non contradiction behind, but I prefer to do it for the sake of less-incompleteness. The non contradiction in systems is just a prejudice. The systems are incomplete, so the contradictions are caused by their own limits. Extending the system allows us to contain the paradoxes, and once we <understand|comprehend> paradoxes, they stop being contradictory.

Perhaps the first reason for entering happily the realm of paradoxes and contradictions is because we know that our systems are incomplete. This means that, even when they might be valid for our present context, as our context enlarges, our systems sooner or later will not be able to be consistent with all our context, as it has happened all through history. And also because we are aware that people who do not share crucial parts of our context, they will not agree with our systems. We are predicting the failure of our systems outside our context (as 1+1=2 in a decimal context, but 1+1=10 in a binary context) by perceiving its limits. This means that we are predicting that someone will say: "Your system is wrong" (related to his incompatible context). But by saying this, he proves that we were right. But people inside our context will say: "You (re)are right", which as our context evolves will prove to be "wrong" (or less-incompletely speaking, "less-incomplete"). This is a paradox. But this paradox makes us less-incomplete, because we are containing the contexts of the people who say "You are wrong". Our theory will re-be true and false at the same time (related to different contexts), but it will **always** re-be true and false at the same time, as compared with theories which now claim to "be" true, without admitting that outside their context they might re-be false.

We should not only seek for the truth or falsity or indeterminacy of an idea, because these are relative (absolute truth (and falsity) is unreachable (but we can approximate as much as we want to)) but for its **less-incompleteness**. This less-incompleteness is also relative, based on our beliefs (one could argue that a small context could be more close to the a-being than another one which contains it, but this cannot be discussed, just believed).

A less-incomplete theory should contain theories which are not wrong, but incomplete in our context (i.e. myths, religions, dogmas, etc...). To be against something is useful, but it is more useful to contain both (<A> and <Anti-A>'s). For example, a less-incomplete political theory should contain socialism and capitalism, despotism and anarchism, dictatorship and democracy, fascism and republic, communism and terrorism... At least, understanding why each one exists and because of what, and for what each one is more suitable; take the things WE need from them, reject the rest, add a bit of sugar, mix it for five minutes (or until the foam is assimilated), and you have another utopia! *If you know the rules of the game, you can change them...*

### 6.1. Neutrosophy and Incompleteness

Neutrosophy, as other systems, embraces the spirit just described: attempting to contain as most contexts as possible, even when they are contradictory, for attempting to be as less-incomplete as possible. But of course Neutrosophy itself predicts its own decay. It is not the *non plus ultra*. It is not finished[2]. But <we cannot|there is no need to> go further now because it fills successfully our contexts.

For example, we could add more values (concepts) to a (Neutrosophic) vector than True, Indeterminate, and False. We just need to define them... (e.g. {T, I, F, T＾I, T＾F, I＾F, T＾I＾F, ~[T|I|F]}) How useful it is this? Someone might ask the same about T, I and F. Why not only T & F, or only F. It depends on the context we are, on the things we need. Logics are just a tool. It depends on what we want to do that we need to make a different tool (of course, there is no "ultimate" tool). Since Neutrosophy is also based in axioms, it is also incomplete.

---

[2]Is there such a thing as a "finished idea"???

<A> is finite, but <Neut-A> is infinite. Therefore, the ideas will evolve infinitely. There are infinite <Anti-A>'s, each one related to a different context (Reference System in Smarandache (1995)). The values/ranges of truth, indeterminacy and falsity (T. I. F.), and are dynamic, relative to a context, and there are infinite contexts, so about any event, there are infinite number of Neutrosophic values/ranges, and any context a-is also infinite (it cannot be **completely** described, as well as the event) Which one is the more representative? The answer to this question is also related to a context! All the "answers" are related to a context, closed, finite... the interesting thing is that all questions seem to be open and infinite... so they can be answered in an infinitude of different (incomplete) ways.

If <A> is combined with <Anti-A>'s, let's suppose it evolves into <Neut-A>. But this is infinite, so it would not be completely <Neut-A>, but as noted by Smarandache (1995), <Neo-A>. Each step/cycle the idea is less-incomplete, but there will always be an infinitude of <Anti-A>'s, and an unreachable <Neut-A>... From a smaller context, <Neo-A> would look just as the old <A>?

## 7.Contextuality

Things re-are in dependence of their context. Since there is an infinitude of contexts, things can re-be in an infinitude of ways. It is only inside a specific context, which we can speak about the truth or falsity or indeterminacy of a proposition.

### 7.1. Context-dependant logic

Every proposition P can only have a truth value (or vectors) *in dependence* of a context C. This truth value/vector is *relative* to the context C. Propositions have only sense (in the sense of Frege (1892)) inside a context(s). Propositions have no sense without a context.

The truth values/vectors of a given proposition can change with context. So, for example "This proposition is false" has a value of 0.5 in Lukaciewicz logic, [1,1] in multidimensional logic, (1,1,1) in Neutrosophic Logic, and "?!" in Aristotelean logic. Or, the proposition "The king of France wears a wig" would be, in terms of multidimensional logic, nor true nor false ([0,0]) in the XX[th] century context, but true ([1,0]) in the context of the 1[st] of January of 1700.

We are just indicating the limits of logics. Logics are just tools. They are useful only inside a context. The context determines the logic. If a proposition goes beyond the context, the logic developed in the context will not be able to contain it (but not necessarily vice versa).

7.2. Derivations of contextuality

Since all <Anti-A>'s of <A> are related to a reference point, we can find all reference points so that all the elements of <Neut-A> are contradictory with <A> (and later with themselves...). This is, any element of <Neut-A> can be potentially <Anti-A>, you just need to have a reference point.

There will always a-be injustice, because this one is relative. Since different people have different contexts (or we can use the word Seelenzustand (soul state), to refer to the personal context, to distinguish from a general context)... So, since people have different

Seelenzustandes, we cannot speak of absolute justice, so things will be just for the people with power... The less-catastrophic panorama (and most naive...) would be that the people in the power would have the less-incomplete Seelenzustandes, trying to contain and understand as many Seelenzustandes as they can, so, if they are just, in spite their relativity, they will be just as well for all the people whose Seelenzustandes they contain.

If ideas are different in each Seelenzustand... well, they might have many similarities, but on the other hand there is the problem of language representing ideas.... one can quote the words of another in another context to communicate different ideas... but the "problem" of language is a different story... ("In Philosophy there are no problems, just opinions (like this one)")

All adjectives are relative. Thus, we can find "opposite" (related to a reference point) adjectives for the same object from different perspectives. (What is wrong from one perspective is right from another, what is beauty-ugliness, good-bad, complex-simple, complicated-simple, complicated-non-complicated, etc...). Are there adjectives which do not behave this way?

## 8. Metacontextuality

Following the ideas exposed above, we can argue that, in order to be less-incomplete, we need to strive for a *metacontextuality*, containing as many contexts as possible. We **believe** this is the only way we have to approach to the a-being more than we already have, but of course there might be other ways, and the proposed way will not be definite.

This does not mean that we will agree with "silly" ideas. This means that instead of just declaring them silly and forgetting about them, we will try to understand what led people to have such ideas, in order to try to see which perspective of the a-being they had. Then our perspective will be greater than if we would just ignore the ideas we do not agree with.

We are tied to a context, a relative one. This is because there re-are no basic elements from where we can build the rest of our world. There are only more complexity and indeterminism in the entrails of the subatomic particles and quarks. There a-is no essence in the universe, because everything is related. Everything a-is for and because of everything. Everything is based on everything. It is naive to try to justify the world once we are already on it. Everything is a condition of everything. Otherwise, it would not be AS IT IS.

Therefore everything can re-be seen in terms of everything (just as a Turing Machine (Turing, 1936) can represent all computations (and computations can represent Turing Machines)). There a-is no base. Everything can re-be a base.

Metacontextuality, as the one Neutrosophy and other currents strive for, by consequence predicts *tolerance*. This is, if we try to contain as much contexts as possible, we will be able to be less intolerant to contexts and Seelenzustandes that are neglected from absolutist non contradictory points of view. And this tolerance should be able to prevent many *conflicts*.

But may metacontextuality lead to *indifference*? This is, since everything can re-be depending on a specific context, does it matter which context we choose? I believe that this might be avoided if we put reason on its place, leaving place for experience and beliefs. Anyway this issue should not be disregarded.

## 9. Conclusions

I **believe** that the question "Do you **believe** in an absolute reality/truth?" is on the same level as "Do you believe in God?". This is, the question is completely metaphysical... And relative truths are incomplete. Should we keep searching for truth of things? I believe now it would be more useful to search for the **less-incompleteness** of things. This can be achieved by enlarging our contexts, containing as many contexts as we can.

Our limits are the ones we set to ourselves. We need only to take the blindfolding of our prejudices in order to attempt a **metacontextuality**. And Neutrosophy bravely does this. If the ideas exposed here, and the ones of Neutrosophy, are not assimilated by our society, this will be because the society does not need them. This is, it "functions" based precisely on the partial blindness of the individuals. Then, should we try to help everyone to open their eyes? We already are doing so, they will open their eyes only if they *want* to.

## 10. Acknowledgements

This work was supported in part by the Consejo Nacional de Ciencia y Tecnologia (CONACYT) of Mexico.